\documentclass[journal]{IEEEtran}
\IEEEoverridecommandlockouts
\usepackage[utf8]{inputenc}
\usepackage[top=2cm, bottom=2cm, left=2cm, right=2cm]{geometry}
\usepackage{float}
\usepackage{caption}
\usepackage{subcaption}
\usepackage{circuitikz}
\usepackage{tikz,pgfplots,pgfplotstable,booktabs}
\usepackage[official]{eurosym}
\usepackage{amsmath,amssymb,amsfonts}
\usepackage[mathscr]{euscript}
\usepackage{xcolor}
\usepackage{colortbl}
\pgfplotsset{compat=1.15}
\usepackage{comment}
\usepackage{multirow,array}
\usepackage{cite}
\usepackage{algpseudocode}
\usepackage{algorithm}
\usepackage{graphicx}
\usepackage{textcomp}
\usepackage{mathtools}

\DeclarePairedDelimiter\floor{\lfloor}{\rfloor}

\newenvironment{ldescription}[1]
  {\begin{list}{}%
   {\renewcommand\makelabel[1]{##1\hfill}%
   \settowidth\labelwidth{\makelabel{#1}}%
   \setlength\leftmargin{\labelwidth}
   \addtolength\leftmargin{\labelsep}}}
  {\end{list}}

\allowdisplaybreaks
\begin{document}
%%%%%%%%%%%%%%%%%%%%%%%%%%%%%%%%%%%%%%%%%%%%%%%%
%\title{Integrating Automatic and Manual Reserves  in Optimal Power Flow via Chance Constraints}
%\title{Limiting the Probability of Adverse Operations: \\ Combining Chance-Constrained and Robust Optimal Power Flow for Dispatch and Reserve Scheduling}
%\title{\color{blue}Limiting the Probability of Adverse Operations: \\ Combining Chance-Constrained and Robust Optimal Power Flow for Reserve Scheduling}
\title{\color{black} Unifying Chance-Constrained and Robust Optimal Power Flow for Resilient Network Operations}
\author{\'Alvaro Porras,  Line Roald, Juan Miguel Morales and Salvador Pineda
\thanks{\'Alvaro Porras, Juan Miguel Morales and Salvador Pineda are with the OASYS research group, University of Malaga, Malaga, Spain. E-mail: alvaroporras19@gmail.com; juanmi82mg@gmail.com (corresponding author); spineda@uma.es. Line Roald is with the Dept. of Electrical and Computer Engineering, University of Wisconsin-Madison, Madison, Wisconsin, USA. E-mail: roald@wisc.edu.}
\thanks{This project has received funding in part by the Spanish Ministry of Science and Innovation (AEI/10.13039/501100011033) through project PID2020-115460GB-I00 and in part by the European Research Council (ERC) under the European Union's Horizon 2020 research and innovation programme (grant agreement No 755705). \'Alvaro Porras is also financially supported by the Spanish Ministry of Science, Innovation and Universities through the university teacher training program with fellowship number FPU19/03053.
Line Roald is supported by the Department of Energy, Office of Science, Office of Advanced Scientific Computing Research, Applied Mathematics program under Contract Number DE-AC02-06CH11347.
The authors thankfully acknowledge the computer resources, technical expertise and assistance provided by the SCBI (Supercomputing and Bioinformatics) center of the University of Malaga.}}
%%%%%%%%%%%%%%%%%%%%%%%%%%%%%%%%%%%%%%%%%%%%%%%%
\maketitle
%%%%%%%%%%%%%%%%%%%%%%%%%%%%%%%%%%%%%%%%%%%%%%%%
\begin{abstract}
Uncertainty in renewable energy generation has the potential to adversely impact the operation of electric \textcolor{black}{networks}.
Numerous approaches to manage this impact have been proposed, ranging from stochastic and chance-constrained programming to robust optimization. However, these approaches either tend to be conservative or leave the system vulnerable to low probability, high impact uncertainty realizations.
To address this issue, we propose a new formulation for stochastic optimal power flow that explicitly distinguishes between “normal operation”, in which automatic generation control (AGC) is sufficient to guarantee system security, and “adverse operation”, in which the system operator is required to take additional actions, e.g., manual reserve deployment.
% This formulation can be understood as a combination of a joint chance-constrained problem, which enforces that AGC should be sufficient with a high probability, with specific treatment of large-disturbance scenarios that allows for additional, more complex redispatch actions.
The new formulation has been compared with the classical ones in a case study on the IEEE-118 and \textcolor{black}{IEEE-300 bus systems}. We observe that our consideration of extreme scenarios enables solutions that are more secure than typical chance-constrained formulations, yet less costly than solutions that guarantee robust feasibility with only AGC.
\end{abstract}

%%%%%%%%%%%%%%%%%%%%%%%%%%%%%%%%%%%%%%%%%%%%%%%%
%%%%%%%%%%%%%%%%%%%%%%%%%%%%%%%%%%%%%%%%%%%%%%%%
\section*{Nomenclature}
The main notation used throughout the text is stated below for  quick  reference.
%%%%%%%%%%%%%%%%
\subsection{Sets}
\begin{ldescription}{$xxx$}
\item [$\mathcal{G}$] Set of generating units, indexed by $g$.
\item [$\mathcal{L}$] Set of transmission lines, indexed by $l$.
\item [$\mathcal{N}$] Set of nodes, indexed by $n$.
\end{ldescription}
%%%%%%%%%%%%%%%%
\subsection{Parameters}
\begin{ldescription}{$xxxxx$}
\item [$c_g$] Linear operating cost of generating unit $g$ [\euro/MWh].
\item [$c^{d}_g/c^{u}_g$] Downward/Upward reserve capacity cost of generating unit $g$ [\euro/MW].
\item [$c^{-}_g/c^{+}_g$] Downward/Upward reserve cost of generating unit $g$ [\euro/MWh].
\item [$B_{ln}$] Power transfer distribution factor (PTDF) of transmission line $l$ with respect to node $n$.
\item [$d_n$] Forecasted demand at node $n$ [MW].
\item [$\overline{f}_l$] Maximum capacity of transmission line $l$ [MW].
\item [$\underline{p}_g/\overline{p}_g$] Minimum/Maximum output of unit $g$ [MW].
\item[$\overline{r}^{d}_g/\overline{r}^{u}_g$] Ability of generator $g$ to provide downward and upward reserves [MW].
\item [$\tilde{w}_n$] Actual wind power production at node $n$ [MW].
\item [$\widehat{w}_n$] Forecasted wind power production at node $n$ [MW].
\item [$\omega_n$] Error of the predicted wind power at node $n$ [MW].
\end{ldescription}
%%%%%%%%%%%%%%%%
\subsection{Variables}
\begin{ldescription}{$xxxxxxxxxx$}
\item [$p_g$] Power output dispatch of unit $g$ [MW].
\item [$r_g(\boldsymbol{\omega})$] Reserve deployed by unit $g$ [MW].
\item [$r^{-}_g(\boldsymbol{\omega})/r^{+}_g(\boldsymbol{\omega})$] Downward/Upward reserve deployed by unit $g$ [MW].
\item [$r^{d}_g/r^{u}_g$] Downward/Upward reserve capacity of unit $g$ [MW].
\item [$\alpha_g(\boldsymbol{\omega})$] Manual adjustment of power output dispatch at unit $g$ [MW].
\item [$\beta_g$] Participation factor of unit $g$.
\end{ldescription}

%%%%%%%%%%%%%%%%%%%%%%%%%%%%%%%%%%%%%%%%%%%%%%%%%%%
\section{Introduction}
\label{sec:introduction}

\IEEEPARstart{O}{ptimization} \textcolor{black}{under uncertainty is a challenging task, and even more so in network systems where uncertainty realizations in different parts of the network can collude to create unexpected difficulties in maintaining balanced network conditions and managing flow constraints. In many systems, it is possible to distinguish between \emph{normal operations}, where the system is operating as desired and (minor) disturbances can be managed using simple, automatic controls, and \emph{adverse operations}, where system security is challenged and additional control actions may be necessary to ``save'' the system from significant impacts. It is typically desirable to limit the probability that the system enters into an adverse operating condition, i.e. to ensure that the normal operational controls are sufficient to guarantee constraint satisfaction with a high probability. However, it is also important to ensure that there exist effective controls to limit system impacts and recover feasible operations if we enter into an adverse operating condition. In this paper, we address this problem by exploring a problem formulation that combines features of chance-constrained and robust optimization, and is related to the idea of chance-constrained programming with recovery proposed in \cite{liu2016decomposition}. %Chance-constrained optimization is used to limit the probability of entering into adverse operations, while robust optimi to ensure system security both during normal and adverse operations.
To make our problem formulation concrete, we consider the problem of operating an electric power grid under uncertainty. However, similar problems of managing the probability of adverse operations and ensuring feasibility recovery in all scenarios could be applied to a range of other networks, including e.g. supply chain networks, emergency response networks, and others.
}
%\IEEEPARstart{T}{he}

\subsection{\color{black}Electric Grid Operation Under Uncertainty}
The
Optimal Power Flow (OPF) is a classical tool widely used for day-ahead and real-time power system operations, electricity markets, long-term planning, and many other applications \cite{shahidehpour2003market}. In its deterministic version, the OPF problem seeks to determine the least-costly dispatch of thermal generating units to satisfy the system's demand, while complying with the technical limits of production and transmission network equipment \cite{frank2012optimal}.
However, the increasing integration of renewable energy sources into power systems leads to increased variability and uncertainty in both the power generation and associated power flows. Understanding and quantifying the impact of this uncertainty on decision-making problems such as the OPF is crucial to ensure the secure operation of power systems \cite{xie2010wind}.  %As a result, operators need to make decisions made in the OPF problem, since there is an increased need for balancing demand through reserves \cite{mohandes2019review}.

%\cite{carpentier1979optimal, stott2012optimal,}

%Given the inherently uncertain nature of the new energy sources dominating power systems, e.g., PV and wind power,
Given this context, a large and growing body of work
has addressed the stochastic version of the optimal power flow (SOPF) problem \cite{roald2023power}. The SOPF aims to minimize expected operational cost and avoid constraint violations while considering the uncertainty in its random parameters. Existing works deal with uncertainty in SOPF using different approaches such as multi-stage stochastic programming \cite{morales2009economic}, robust or worst-case optimization \cite{jabr2015robust,warrington2013policy, lorca2018adaptive} or chance-constraints  \cite{vanackooij2011, vrakopoulou2013, bienstock2014chance, lubin2015robust, hou2020tunning, esteban2023distributionally}. The major challenge is to design a model that captures the risk of constraint violations and accurately reflects the operation of power systems, while maintaining computational tractability.

\subsection{\color{black} Activation of Generation Reserves}
When modeling the impact of uncertain generation on short-term operations (i.e. day-ahead until real-time), it is common to assume that forecast errors and renewable energy variability will be balanced by the deployment of generation reserves, and in particular, by systems such as the \emph{automatic generation control} (AGC) \cite{warrington2012robust}.
A benefit of modeling system balancing through AGC is that it is naturally represented as an affine control policy, which also simplifies the solution of the optimization problem.
% \cite{vrakopoulou2013,roald2013analytical,lubin2015robust,pena2021jointopf}
%As a common approach of most of the existing papers working on SCOPF is to model the generation adjustments, which cope fluctuations on demand or renewable generation, i.e., reserve deployment, as an affine control policy \cite{warrington2012robust}. As a consequence, two advantages are achieved, both modeling the use of \emph{automatic generation control} (AGC) for short-term system balancing and resulting in tractable optimization problems.
%
However, the common AGC models implemented in the literature typically %present several drawbacks. First, they
assume that all generators contribute reserve power according to the affine policy, even for large uncertainty deviations. In reality, generator output will saturate (or stop increasing/decreasing as the deviation grows larger) when they hit their lower or upper generation limit. Furthermore, operators generally take additional actions to manage both balancing and congestion in situations with very large deviations. For example,  the \emph{North American Reliability Corporation} (NERC) \cite{nerc} standard for regulating the use of AGC, BAL-005, states that if the AGC becomes inoperative or may impair the reliability of the interconnection, the system operator must use manual controls to adjust generation in order to guarantee balance.

While a limited number of studies have shown that modeling generation saturation \cite{kannan2020stochastic} or accounting for manual reserve activation during large deviations \cite{roald2015optimal} leads to better operating conditions, these models are often computationally expensive.
%(resulting in high operating costs), whereas, in real power systems, generators have the capability to saturate their power output dispatch as it is stated in \cite{kannan2020stochastic}. The latter, the affine AGC models fail to model the actions of system operators during more severe operating conditions, where manual adjustments/actions are typically performed by system operators to restore system balance or to relieve congestion and overloads.
% Specifically, \emph{North American Reliability Corporation} (NERC) \cite{nerc} developed a standard for regulating the use of AGC, BAL-005, which states that if the AGC becomes inoperative or may impair the reliability of the interconnection, the system operator must use manual controls to adjust generation in order to guarantee balance.
%
%As a countermeasure to the disadvantage of modeling the affine control policy for all operating conditions, it is possible to
Thus, a more common approach is to apply the affine control policy, but
explicitly disregard constraint satisfaction in a fraction of the most severe operating conditions. This is typically done by introducing chance constraints that allow violations in a (typically small) percentage of scenarios \cite{vanackooij2011, vrakopoulou2013, bienstock2014chance, lubin2015robust, hou2020tunning, esteban2023distributionally}, or by solving robust optimization formulations where the uncertainty set has been designed to contain a certain probability mass \cite{Kostas2013}.
%\cite{miller1965chance, vrakopoulou2013, lubin2015robust, hou2020tunning, Vrakopoulou2013acopf,vrakopoulou2013reserve}. %Examples from the literature on \emph{joint chance-constrained} OPF (JCC-OPF) \cite{vanackooij2011}, given its benefits, there is a vast related literature \cite{bienstock2014chance,lubin2015robust,hou2020tunning,pena2021jointopf,esteban2023distributionally} where the affine control policy, AGC, is applied to manage uncertainty fluctuations.
%
% The probabilistic (chance-constrained) version of SOPF problem is commonly modeled in the literature by means of a joint chance-constraint {\color{red} [references]}, where the constraint satisfaction is simultaneously disregarded for generators and transmission lines constraints. This modeling choice seems more intuitive than a single chance constraint per element which implies a special expertise of each element of the system in order to determine its probability of violation.

Unfortunately, by failing to model the impact of the worst scenarios (those for which the constraint satisfaction is discarded), a chance-constrained formulation may leave the system vulnerable to large disruptions that include generator and line outages, or load shed. As discussed in \cite{bienstock2014chance}, there could be instances where the combination of generators and renewable outputs collude to produce power flows that significantly exceed the nominal line ratings, even in the absence of a large total power deviation.
When the maximum rating of a line is exceeded, this line becomes more likely to trip, leaving the network vulnerable to cascading failures and associated load shed.

%In the literature, several works try to deal with these extreme events by mitigating their occurrence through robust optimization as \cite{jabr2015robust,lorca2018adaptive}. Note that the decisions in robust approaches are optimal for the worst case what is too conservative.

%Authors in \cite{roald2015optimal} model a manual adjustment in the generation for power imbalances exceeding a certain threshold, i.e., a huge imbalance. Otherwise, the power imbalance is corrected by performing the AGC. However, there are higher risky or costly, in a nutshell, adverse scenarios that require this type of action, i.e., manual adjustments, without having a large imbalance. For example, given the spatial distribution of the sources of uncertainty, a large congestion event can occur in the system with a small total imbalance. Roald \emph{et al.} \cite{roald2016corrective} design corrective actions to handle post-contingency events which can be implemented through automatic (AGC) or manual adjustments. However, this research rules out the possibility of the AGC implementing an unreliable action in the face of a scenario of high variability or uncertainty in renewable plants. This would imply the use of a manual adjustment by the system operator without the need to be a contingency in the system.

\subsection{\color{black}Contributions}
To address this issue, in this work we \textcolor{black}{propose a new SOPF formulation that distinguishes between two different operating regimes, namely \emph{normal operation} and \emph{adverse operation}. In the SOPF context, \emph{normal operation} refers to a situation in} which AGC is sufficient to maintain the system balance, while \emph{adverse operation} \textcolor{black}{refers to a situation in which} the system operator may need to implement additional actions, such as manual adjustments, to preserve system security.
%To address this issue, in this work we propose a new SOPF formulation where, in the worst-case scenarios, operators can leverage additional, manually activated reserves to reduce their impact. The proposed formulation distinguishes between two different operating regimes: \emph{normal operation} and \emph{adverse operation}. In \emph{normal operation}, AGC is sufficient to maintain the system balance, while in \emph{adverse operation}, the system operator may need to implement additional actions, such as manual adjustments, to preserve system security.

Unlike the standard \emph{joint chance-constrained} OPF (JCC-OPF), which limits the joint probability of violation of technical constraints, our formulation uses a joint chance-constraint to \emph{control} the probability of utilizing different reserve actions. Thus, we can impose that AGC alone is to be sufficient with a high probability, while additional corrective actions are only to be implemented for the most adverse scenarios. In doing so, we reduce the need for frequent manual intervention by operators (computational expensive), while also guaranteeing that additional resources are available to handle adverse operating conditions, e.g., by scheduling more reserve capacity for manual deployment.

To demonstrate the suitability of our proposed formulation, we conduct a computational experiment that compared it to two state-of-the-art approaches. The former is the standard JCC-OPF, whose drawback is to leave the power system vulnerable to severe events, e.g., by dispatching insufficient generation capacity or giving up on alleviating congestion. The second approach guarantees robust feasibility using AGC only, which results in conservative solutions with increased operating costs, for instance, due to inefficient and oversized generation capacity. Our novel formulation results in decisions that are more reliable than the former approach and more cost-efficient than the latter.

The main contributions of our work are thus the following:

\begin{itemize}
\item We propose a novel optimal power flow formulation that accounts for \textcolor{black}{both adverse and normal operations and their relation to the} various reserve actions employed in actual power system operations, such as AGC and manual re-dispatch.
\item \textcolor{black}{Our problem formulation combines aspects of chance-constrained and robust optimization.} We use a joint chance constraint to restrict the probability \textcolor{black}{of entering adverse operating conditions, while enforcing that feasibility recovery is possible for all scenarios with the use of additional controls. Specifically,} we limit the probability of manual adjustments occurring instead of limiting the probability of violation of technical constraints, \textcolor{black}{and guarantee that the use of manual controls ensures that the system remains feasible in all scenarios. This is representative of} %Hence,
a realistic setting where the AGC operates under ordinary system conditions and the manual adjustments, which are not automatic, are implemented under adverse scenarios only.
\item We show that our approach \textcolor{black}{provides an opportunity to obtain solutions that are different than existing chance-constrained and robust optimization approaches. In particular, our approach} yields solutions that are more reliable than the conventional joint chance-constrained DC-OPF, yet less costly than those approaches that guarantee robust feasibility \textcolor{black}{with AGC alone}.
\end{itemize}

The rest of this paper is organized as follows. Section~\ref{sec:formulation} describes the proposed SOPF formulation, which is derived from two traditional approaches in the literature. Section~\ref{sec:solution} describes the reformulation and algorithm used to make our proposal tractable and computationally efficient. Section~\ref{sec:eval} explains the methodology we use to benchmark our approach, while Section \ref{sec:case} discusses experimental results from a case study. Finally, conclusions are drawn in Section~\ref{sec:conclusion}.

%%%%%%%%%%%%%%%%%%%%%%%%%%%%%%%%%%%%%%%%%%%%%%%%
\section{Problem Formulation}
\label{sec:formulation}
We start this section by introducing a standard and well-known formulation of the \emph{joint chance-constrained DC-OPF problem}. \textcolor{black}{We use this formulation a basis to construct and  motivate our proposed stochastic OPF formulation, which is presented immediately after. }

%Before diving into the mathematical formulations, which is formulated to capture the operational characteristics of electric power systems, we give a brief overview of the rationale behind the proposed formulation and place it in context of more general stochastic optimization.
\textcolor{black}{Before diving into the detailed mathematical formulations, we would like to highlight some distinguishing features between our proposed formulation and common stochastic optimization approaches. While our proposed formulation contains a joint chance constraint and thus is a chance-constrained program, there are some significant differences to common joint chance-constrained optimization. In our formulation, the chance constraint controls the probability of manual reserve activation, not the probability of constraint violation. In a more general setting, the chance constraint can be seen as managing the probability of switching from normal to adverse operations. Furthermore, our proposed problem formulation guarantees that all scenarios can be made feasible by leveraging the additional controls (in our case manual reserve activation), and thus provides robust constraint satisfaction across scenarios. As such, our proposed formulation combines features of chance-constrained and robust optimization.}

%\subsection{Notation}
{\color{black}For notational simplicity, as in \cite{hou2020tunning}, we assume that there is one dispatchable generator, one uncertain power source (e.g., a wind farm) and one power load per node $n$}. The power dispatch of the generator, the power produced by the uncertain power source and the power consumed by the power load are denoted by $p_n$, $\tilde{w}_n$, and $d_n$, respectively. %\textcolor{red}{1-hour period, MW = MWh}

The power $\tilde{w}_n$ generated by the uncertain power source at node $n$ is modeled as a random variable, which we decompose as  $\tilde{w}_n = \widehat{w}_n + \omega_{n}$, with $\widehat{w}_n$ being the forecast value (assumed unbiased) and $\omega_{n}$ the associated forecast error. The system-wide aggregate forecast error is given by $\Omega = \sum_{n \in \mathcal{N}} \omega_n$, \textcolor{black}{and translates into a system-wide power imbalance} which is balanced by the dispatchable generators through the deployment of reserve. The reserve deployment follows an affine control policy, modeling the actions of the \emph{Automatic Generation Control} (AGC). According to this policy, the reserve provided by the generator at node $n$ is given by
$r_n(\boldsymbol{\omega}) = -\beta_n\Omega$, where $\beta_n$ is the participation factor and $\boldsymbol{\omega}$ represents the vector of power forecast errors across all nodes. Furthermore, we distinguish between \emph{upward} or \emph{positive} reserve $r^+_n(\boldsymbol{\omega})$ and \emph{downward} or \emph{negative} reserve $r^-_n(\boldsymbol{\omega})$, with $r_n(\boldsymbol{\omega}) = r^+_n(\boldsymbol{\omega}) - r^-_n(\boldsymbol{\omega})$.

%\subsection{Preliminaries: Joint Chance-Constrained OPF}
With this notation in place, the \emph{joint chance-constrained DC-OPF problem} can be formulated as follows:
\begin{subequations}
\label{eq:JCC-OPF}
\begin{align}
%Objective Function
\min_{\Xi}  &  \quad\sum_{n} c_{n} \, p_{n} + c^{u}_{n} r^{u}_{n} + c^{d}_{n} r^{d}_{n} \notag\\
& \hspace{2cm} + \mathbb{E} \left[c^{+}_n r^{+}_n(\boldsymbol{\omega}) - c^{-}_n r^{-}_n(\boldsymbol{\omega})
 \right] \label{eq:JCCOPF_objective}\\
% Power Balance
\text{s.t.} & \sum_{n} \beta_n = 1 \label{eq:JCCOPF_balance1}\\
& \sum_{n} (p_n + \widehat{w}_n - d_{n}) = 0 \label{eq:JCCOPF_balance2}\\
& \underline{p}_{n} + r^{d}_{n} \leq p_n  \leq \overline{p}_{n} - r^{u}_{n}, \quad \forall n  \label{eq:JCCOPF_gen}\\
&  0 \leq r^{d}_{n} \leq \overline{r}^{d}_n, \quad \forall n  \label{eq:JCCOPF_maxdores}\\
&  0 \leq r^{u}_{n} \leq \overline{r}^{u}_n, \quad \forall n  \label{eq:JCCOPF_maxupres}\\
& r^{+}_n(\boldsymbol{\omega}) - r^{-}_n(\boldsymbol{\omega}) = -\Omega\beta_n, \forall n \label{eq:JCCOPF_updownreserve}\\
% Joint Chance Constraint
& \mathbb{P}
\left(\begin{array}{l}
    -r^{d}_{n} \leq  -\Omega\beta_n \leq r^{u}_{n}, \quad \forall n  \\
    -\overline{f}_{l} \leq \displaystyle\sum_{n} B_{ln} (p_{n} -\Omega\beta_n +\\
    \hspace{0.7cm}+ \widehat{w}_n + \omega_n - d_n ) \leq \overline{f}_{l}, \forall l
\end{array} \right) \geq 1 - \epsilon \label{eq:JCCOPF_jointCC}\\
& \beta_n, r^{+}_n(\boldsymbol{\omega}), r^{-}_n(\boldsymbol{\omega}) \geq 0, \quad \forall n \label{eq:JCCOPF_charac},
\end{align}
\end{subequations}
\noindent where $\Xi = (p_n, r^{+}_n(\boldsymbol{\omega}), r^{-}_n(\boldsymbol{\omega}), r^{d}_{n}, r^{u}_{n}, \beta_n)$ is the set of decision variables. We remark that $r^{d}_{n}$ and $r^{u}_{n}$ are the downward and upward reserve \emph{capacity} provided by the dispatchable generator at node $n$. This reserve capacity is procured by the system operator before the forecast errors $\boldsymbol{\omega}$ are known. %and as such, is independent of these.

The three terms of the objective function \eqref{eq:JCCOPF_objective} to be minimized correspond to the power dispatch cost, the cost of procuring reserve capacity, and the expected cost related to the actual deployment of that capacity, respectively. %, in that order.
The power balance in the system  is guaranteed for all possible realizations of $\boldsymbol{\omega}$ through  equations \eqref{eq:JCCOPF_balance1} and \eqref{eq:JCCOPF_balance2}.
Note that by requiring $\beta_n\geq 0$ in \eqref{eq:JCCOPF_charac}, we enforce that all reserve deployment only acts to counterbalance forecast errors, rather than also allowing redispatch among generators to counter congestion.
Constraints \eqref{eq:JCCOPF_gen} ensure that the power produced \textcolor{black}{and the reserve capacity offered} by dispatchable generators is within their minimum and maximum power output limits $\underline{p}_{n}$ and $\overline{p}_{n}$. \textcolor{black}{This condition is expressed by adjusting the generation limits for the dispatchable generation by the respective reserve capacities each generator provides}.
Constraints \eqref{eq:JCCOPF_maxdores} and \eqref{eq:JCCOPF_maxupres} set a limit on the maximum reserve capacity each generator is willing or able to provide. Equation \eqref{eq:JCCOPF_updownreserve} models the affine control policy for reserve deployment (AGC) we discussed above. Expression \eqref{eq:JCCOPF_jointCC} constitutes the joint chance-constraint system by which the system operator states that the reserves deployed and the line flows must be within their  bounds with a probability greater than or equal to $1-\epsilon$. Accordingly, the parameter $\epsilon$ is the maximum allowed probability of constraint violation set by the operator \textcolor{black}{to reflect their risk preferences. Note that by implementing a joint chance constraint (as opposed to e.g. multiple single chance constraints), we can interpret $\epsilon$ as a metric for overall system security. In the literature on chance constrained OPF, $\epsilon$ is often chosen to be in the range from 0-10\%.}
In~\eqref{eq:JCCOPF_jointCC}, $l$ is the index of transmission lines in the power network and $\overline{f}_l$ stands for the capacity limit of line $l$.
\textcolor{black}{Note that the joint chance constraint limits the probability that the reserve activation exceeds the contracted reserve capacities, and thus ensures, in combination with \eqref{eq:JCCOPF_gen}, that the generation limits will be satisfied with a high probability.}
Finally,  \eqref{eq:JCCOPF_charac} imposes the positive character of decision variables $\beta_{n}$ and random functions $r^{+}_n(\boldsymbol{\omega})$ and $r^{-}_n(\boldsymbol{\omega})$ for all $n$. Note that the probability in \eqref{eq:JCCOPF_jointCC} is computed over the probability space of $\boldsymbol{\omega}$ and that the equality \eqref{eq:JCCOPF_updownreserve} and the inequality \eqref{eq:JCCOPF_charac} must hold for almost all $\boldsymbol{\omega}$.

The popularity of the joint chance-constrained formulation~\eqref{eq:JCC-OPF} stems from its ability \textcolor{black}{to guarantee overall system security by ensuring that all constraints will remain satisfied with a high probability, while at the same time} reducing the expected system operating cost substantially \textcolor{black}{compared with robust optimization} by allowing the violation of reserve capacity constraints and/or line flow limits under a small $\epsilon$-percentage of realizations of $\boldsymbol{\omega}$. These realizations, or scenarios, are thus the most detrimental to the system in terms of cost. Parameter $\epsilon$ in~\eqref{eq:JCC-OPF} controls the level of risk aversion of the system operator (a lower $\epsilon$ implies more risk averse). If $\epsilon$ is set to $0$, formulation~\eqref{eq:JCC-OPF} becomes \emph{robustly feasible}, meaning that all the constraints and variable limits are to be satisfied with probability one.

%\subsection{Robust Feasibility}
\textcolor{black}{While chance-constrained OPF has gained widespread attention and is closely tied to existing criteria for reliability and reserve procurement}, the critical nature of power systems practically forces operators to guarantee robust feasibility.
%In this regard, formulation~\eqref{eq:JCC-OPF}, even if $\epsilon = 0$, offers an incomplete picture of how power systems are actually operated.
Indeed, in those very few $\boldsymbol{\omega}$-scenarios for which AGC is unable or too costly to ensure the system's integrity, the operators can still take over the affine control policy and \emph{manually} set new operating points for some generators in the system, those needed to guarantee the satisfaction of the system's constraints ideally at the minimum cost. The fact that formulation~\eqref{eq:JCC-OPF} ignores the possibile need for a manual control taking over AGC (albeit with a low occurrence) causes it to underestimate the operating cost when $\epsilon > 0$ or overestimate it when robust feasibility is pursued ($\epsilon = 0$). The ultimate result is that formulation~\eqref{eq:JCC-OPF} may produce uneconomical or suboptimal affine control policies.

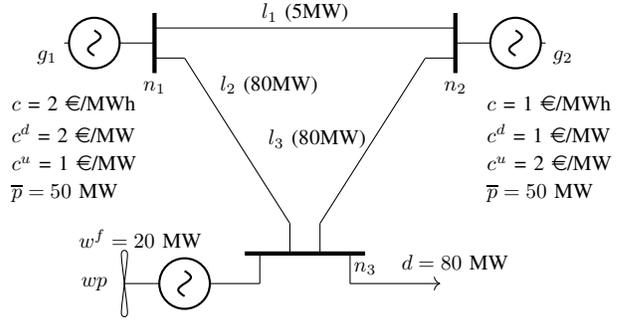
\begin{figure}
\centering
\begin{circuitikz}[scale=0.8,every node/.style={scale=0.8}]
%Nodes
\draw [ultra thick] (0,4)  -- (0,3) node[anchor=north]{$n_{1}$};
\draw [ultra thick] (5,4)  -- (5,3) node[anchor=north]{$n_{2}$};
\draw [ultra thick] (1.5,0)  -- (3.5,0) node[anchor=north]{$n_{3}$};
%Generators
\draw[] (-1.5,3.5) node[anchor=north east]{$g_1$} to [sV] (-0.5,3.5);
\draw (-0.5,3.5) -- (0,3.5);
\draw[] (6.5,3.5) node[anchor=north west]{$g_2$} to [sV] (5.5,3.5);
\draw (5.5,3.5) -- (5,3.5);
%windfarm1
\draw[] (0,-0.5)  to [sV] (1,-0.5);
\draw (1,-0.5)--(1.75,-0.5);
\draw (1.75,-0.5)--(1.75,0);
\draw plot [smooth cycle] coordinates {(-0.55,-1) (-0.5,-0.5) (-0.45,0) (-0.55,0) (-0.5,-0.5) (-0.45,-1)};
\draw (-0.5,-0.5)--(0,-0.5);
\node at (-1,-0.5) {$wp$};
%Lines
%l1
\draw (0,3.75) -- (5,3.75);
\node[] at (2.5,4) {$l_1$ (5MW)};
%l2
\draw (0,3.25) -- (0.5,3.25);
\draw (0.5,3.25) -- (2.25,0.5);
\draw (2.25,0.5) -- (2.25,0);
\node[] at (1.9,2.8) {$l_2$ (80MW)};
%l3
\draw (2.75,0) -- (2.75,0.5);
\draw (2.75,0.5) -- (4.5,3.25);
\draw (4.5,3.25) -- (5,3.25);
\node[] at (2.7,1.9) {$l_3$ (80MW)};
%Demands
\draw (3.25,-0) -- (3.25,-0.5);
\draw[->] (3.25,-0.5) -- (4.75,-0.5);
\node[anchor=east] at (6.0,-0.15) {$d = 80$ MW};
% Wind farm
\node[anchor=east] at (0.9,0.25) {$w^{f} = 20$ MW};
%Data generators
\node[anchor=west] at (-2.5,2.5) {$c$ = 2 \euro/MWh};
\node[anchor=west] at (-2.5,2) {$c^d$ = 2 \euro/MW};
\node[anchor=west] at (-2.5,1.5) {$c^u$ = 1 \euro/MW};
\node[anchor=west] at (-2.5,1) {$\overline{p} = 50$ MW};
\node[anchor=west] at (5.4,2.5) {$c$ = 1 \euro/MWh};
\node[anchor=west] at (5.4,2) {$c^d$ = 1 \euro/MW};
\node[anchor=west] at (5.4,1.5) {$c^u$ = 2 \euro/MW};
\node[anchor=west] at (5.4,1) {$\overline{p} = 50$ MW};
\end{circuitikz}
\caption{Three-node illustrative example}
\label{fig:3bus}
\vspace{-0.5cm}
\end{figure}

To illustrate our point, we use an example based on the small power system depicted in Fig. \ref{fig:3bus}. The system includes two thermal generating units with the linear production costs, reserve capacity costs and maximum power limits indicated in the figure. For simplicity, the susceptances of all lines are assumed to be 1 p.u. and the capacity of each line is also specified in the figure. A single demand of 80 MW is located at node $n_3$, where there is also a wind farm with a predicted output of  20 MW. We assume that the associated (random) forecast error can take on three different values only, namely, $20$, $10$, and $-20$ MW, corresponding to three equally probable realizations or scenarios 1, 2 and 3, in that order. The costs of deploying upward and downward reserve, i.e., $c^{+}_n$ and $c^{-}_n$, are 1.2 and 0.8 times the generator's linear operating cost, respectively.

\begin{table*}[h]
\caption{Results -- Illustrative Example}
\vspace{-0.2cm}
\begin{center}
\begin{tabular}{lccccccccc}
\hline
Method &$p_1$ &$p_2$ &$\beta_1$ &$\beta_2$ &$r^{u}_1$ &$r^{u}_2$ &$r^{d}_1$ &$r^{d}_2$ & Cost (\euro) \\
\hline
model \eqref{eq:JCC-OPF} ($\epsilon=1/3$)  &15  &45    &0   &1  &0 &0 &0 &20 &95\\
model \eqref{eq:JCC-OPF} ($\epsilon=0$)  &12.5  &47.5    &0.625   &0.375  &12.5 &7.5 &12.5 &7.5 &137.5\\
model \eqref{eq:FJCC-OPF} ($\epsilon=1/3$)  &15  &45    &0   &1  &20 &0 &0 &20 &123\\
\hline
\end{tabular}
\label{tab:Results-IE}
\end{center}
\end{table*}

Results from problem~\eqref{eq:JCC-OPF} when $\epsilon = 1/3$ and $\epsilon = 0$ are shown in the first two rows of Table~\ref{tab:Results-IE}. These results include  the optimal power dispatch, participation factors and procured reserve capacities, together with the optimal expected operating cost. Unsurprisingly, the results are quite sensitive to $\epsilon$. For example, when this is set to zero (to achieve robust feasibility), much more reserve capacity is to be procured than when we allow the system's constraints to be violated under one of the scenarios, in particular, scenario~3. Accordingly, the cost increases from \euro95 to \euro137.5, when $\epsilon$ goes from 1/3 to 0. In contrast, if we take the solution delivered for $\epsilon = 1/3$ and scenario 3 actually occurs, meaning that the wind power forecast error takes on the value $-20$MW, the AGC requires generator $g_2$ to increase its production up to $65$MW, that is, above its maximum output limit, while exceeding the maximum capacity of line $l_1$ too. But what is more important from a practical point of view is that, under such a scenario, the solution to problem~\eqref{eq:JCC-OPF} when $\epsilon = 1/3$ does \emph{not} allow for any \emph{manual re-dispatch} that can restore system feasibility  immediately after, because no upward reserve capacity has been procured beforehand. Consequently, formulation~\eqref{eq:JCC-OPF} may be too risky or too costly.

\begin{figure}[tbp]
\centerline{\includegraphics[scale=0.35]{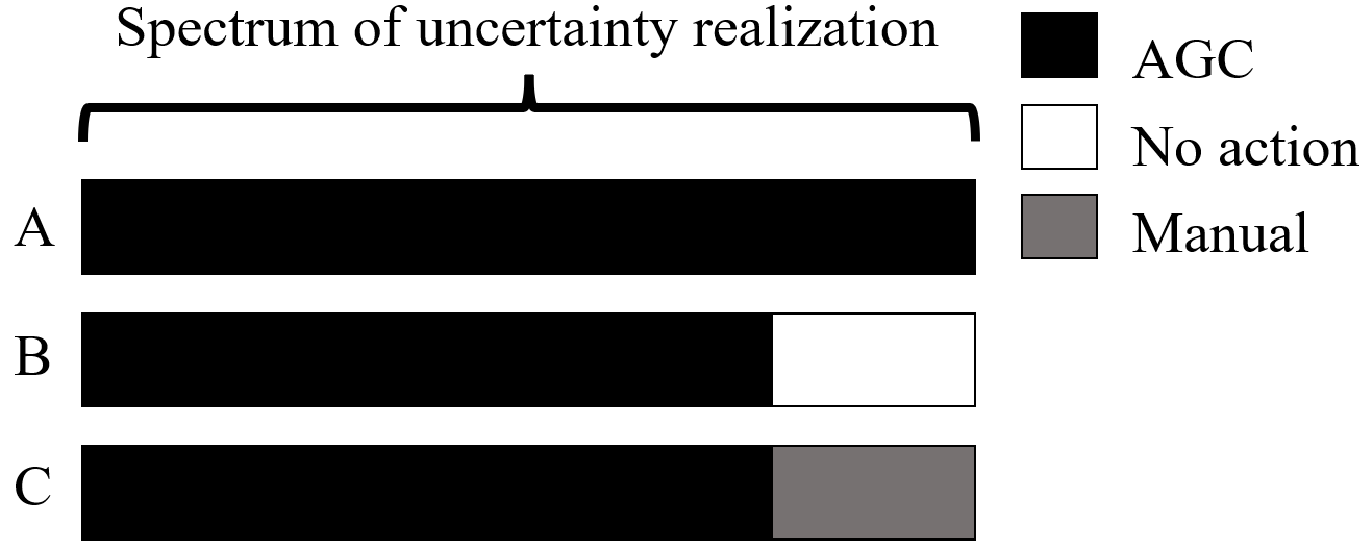}}
\caption{{\color{black}Actions planned over the spectrum of uncertainty realizations to mitigate imbalances and ensure the reliability of the power system.} A: model \eqref{eq:JCC-OPF} with $\epsilon = 0$, B: model \eqref{eq:JCC-OPF} with $\epsilon > 0$, C: Proposal.}
\label{fig:approaches_sopf}
\vspace{-0.5cm}
\end{figure}

To adress this issue, we propose a novel joint chance-constrained formulation of the DC-OPF problem that does account for the possibility of resorting to a manual re-dispatch in this $\epsilon$-percentage of events in which the implementation of AGC is too expensive or even infeasible.
%{\color{black} To provide a more technical description of our proposal, Fig.~\ref{fig:approaches_sopf} illustrates how the mentioned approaches respond to an uncertain scenario.
Our approach involves the use of manual adjustments to balance and ensure reliability under extreme conditions, while model \eqref{eq:JCC-OPF} with $\epsilon=0$ only employs AGC and model \eqref{eq:JCC-OPF} with $\epsilon>0$ does not consider any corrective measures for adverse scenarios. To obtain our formulation, we replace the set of constraints \eqref{eq:JCCOPF_updownreserve}--\eqref{eq:JCCOPF_jointCC} in~\eqref{eq:JCC-OPF} with the following ones:
\begin{subequations}
\label{eq:manual}
\begin{align}
& \sum_n \alpha_n(\boldsymbol{\omega}) = 0 \label{eq:balance_FJCCOPF}\\
& r^{+}_n(\boldsymbol{\omega}) - r^{-}_n(\boldsymbol{\omega}) = -\Omega\beta_n + \alpha_n(\boldsymbol{\omega}), \forall n \label{eq:updownreserve_FJCCOPF}\\
& -r^{d}_{n} \leq  -\Omega\beta_n + \alpha_n(\boldsymbol{\omega}) \leq r^{u}_{n},
\quad \forall n  \label{eq:resca_FJCCOPF}\\
& -\overline{f}_{l} \leq \sum_{n} B_{ln} (p_{n} -\Omega\beta_n + \alpha_n(\boldsymbol{\omega}) +\notag\\
& \hspace{2.5cm}+ \widehat{w}_n + \omega_n - d_n ) \leq \overline{f}_{l}, \quad \forall l \label{eq:flow_FJCCOPF}\\
& \mathbb{P} \left(\alpha_n(\boldsymbol{\omega}) = 0, \quad \forall n\right) \geq 1 - \epsilon  \label{eq:manual_FJCCOPF},
\end{align}
\end{subequations}
\noindent where $\alpha_n(\omega)$ is a random variable that  represents the manual adjustment of the power output of the generator located at node~$n$.

Equation \eqref{eq:balance_FJCCOPF} is required for the implemented power adjustments to preserve the power balance. Equation \eqref{eq:updownreserve_FJCCOPF} is analogous to \eqref{eq:JCCOPF_updownreserve}, but including the manual adjustment requested by the system operator.
Inequalities \eqref{eq:resca_FJCCOPF} and \eqref{eq:flow_FJCCOPF} enforce that the use of AGC in combination with manual re-dispatch guarantees that the system's constraints are satisfied under any realization of $\boldsymbol{\omega}$. Finally, in our proposal, chance-constrained programming is employed for a different purpose than that in~\eqref{eq:JCC-OPF}. Specifically, the chance-constraint \eqref{eq:manual_FJCCOPF} seeks to characterize the use of manual control as an \emph{occassional} recourse action, thus ensuring that AGC is the standard control policy. Again constraints \eqref{eq:balance_FJCCOPF}--\eqref{eq:flow_FJCCOPF} are to be satisfied for almost all $\boldsymbol{\omega}$.

Our proposal can thus be formulated as follows:
\begin{subequations}
\label{eq:FJCC-OPF}
\begin{align}
\min_{\Xi}  \quad &  \sum_{n} c_{n} \, p_{n} + c^{u}_{n} r^{u}_{n} + c^{d}_{n} r^{d}_{n} \notag\\
& \hspace{2cm} + \mathbb{E} \left[c^{+}_n r^{+}_n(\boldsymbol{\omega}) - c^{-}_n r^{-}_n(\boldsymbol{\omega})
 \right] \\
\text{s.t.} \quad & \eqref{eq:JCCOPF_balance1} - \eqref{eq:JCCOPF_maxupres}, \eqref{eq:balance_FJCCOPF} - \eqref{eq:manual_FJCCOPF}, \eqref{eq:JCCOPF_charac}.
\end{align}
\end{subequations}

Coming back to our example, results from~\eqref{eq:FJCC-OPF} are also included in the last row of Table \ref{tab:Results-IE}. Observe that the system operating cost is significantly reduced with respect to that of formulation~\eqref{eq:JCC-OPF} with $\epsilon = 0$. Furthermore, as in the case of formulation~\eqref{eq:JCC-OPF} with $\epsilon = 1/3$, our proposal also renders a solution for which, if scenario 3 occurs, the implementation of AGC violates the maximum output limit of generator $g_2$ and the capacity of line $l_1$. However, unlike the solution to~\eqref{eq:JCC-OPF}, the one delivered by our proposal procures $20$ MW of upward reserve capacity from generator $g_1$ so that the system operator can release this generator from AGC and manually dispatch it at $20$ MW under scenario 3.
\textcolor{black}{For illustrative purposes, Fig.~\ref{fig:approaches_sopf} summarizes how the mentioned approaches respond to an uncertain scenario. }

In the following section, we discuss how we solve formulations~\eqref{eq:JCC-OPF} and~\eqref{eq:FJCC-OPF}.

\section{Solution methodology}
\label{sec:solution}
Chance-constrained programs such as \eqref{eq:JCC-OPF} and \eqref{eq:FJCC-OPF} belong to the class of NP-hard problems. In general,  there is no finite tractable reformulation of the chance-constraint \eqref{eq:JCCOPF_jointCC} or \eqref{eq:manual_FJCCOPF}. As a result, a wide variety of different approaches have been proposed to approximate the feasible region determined by these constraints, namely,  Distributionally Robust Optimization (DRO) \cite{hanasusanto2017ambiguous}, the Scenario Approach (SA) \cite{nemirovski2006scenario}, Sample Average Approximation (SAA) \cite{luedtke2008sample}, and the inner convex approximations based on the Conditional Value-at-Risk (CVaR) \cite{nemirovski2007convex} or ALSO-X \cite{jiang2022also}. In this work, we resort to SAA boosted with bounding, tightening and valid inequalities. Consequently, the chance-constrained programs \eqref{eq:JCC-OPF} and \eqref{eq:FJCC-OPF} are reformulated as mixed-integer programs (MIP). To do that, we assume that  $\boldsymbol{\omega}$ has a finite discrete support defined by a collection of atoms $\{\boldsymbol{\omega}_s \in \mathbb{R}^{|\mathcal{N}|}, s \in \mathcal{S}\}$ and respective probability masses $\mathbb{P}(\boldsymbol{\omega} = \boldsymbol{\omega}_s)=\frac{1}{|\mathcal{S}|}$, $\forall s\in \mathcal{S}=\{1,\dots,|\mathcal{S}|\}$. Accordingly, $\omega_{ns}$ and $\Omega_s$ are realizations of the respective random variables under scenario $s$, and the decisions $\alpha_{ns}$, $r_{ns}^{+}$ and $r_{ns}^{-}$ for the dispatchable unit at node $n$ may vary for each scenario $s$. We define $q = \floor*{\epsilon |\mathcal{S}|}$ and the vector $\boldsymbol{y}$ of binary variables $y_s$, $\forall s\in \mathcal{S}$. Thus, the MIP reformulation of problem \eqref{eq:FJCC-OPF} writes as follows:
\begin{subequations}
\label{eq:SAA_FJCC-OPF}
\begin{align}
\min_{\Theta}  \quad &  \sum_{n} c_{n} \, p_{n} + c^{u}_{n} r^{u}_{n} + c^{d}_{n} r^{d}_{n}  + \frac{1}{|\mathcal{S}|}\sum_s c^{+}_n r^{+}_{ns} - c^{-}_n r^{-}_{ns}
 \\
\text{s.t.}  \quad & \eqref{eq:JCCOPF_balance1}-\eqref{eq:JCCOPF_maxupres} \label{eq:SAA_deterministic}\\
& \sum_{n} \alpha_{ns} = 0, \quad \forall s \label{eq:SAA_balance}\\
& r^{+}_{ns} - r^{-}_{ns} = -\Omega_s\beta_n + \alpha_{ns}, \quad \forall n,s \label{eq:SAA_updownreserve}\\
& -r^{d}_{n} \leq  -\Omega_s\beta_n + \alpha_{ns} \leq r^{u}_{n}, \quad \forall n,s
 \label{eq:SAA_resca}\\
& -\overline{f}_{l} \leq \sum_{n} B_{ln} (p_{n} -\Omega_s\beta_n + \alpha_{ns} +\notag\\
& \hspace{2cm} +\widehat{w}_n + \omega_{ns} - d_n ) \leq \overline{f}_{l}, \quad \forall l,s\label{eq:SAA_flow}\\
& -y_s \overline{r}^{d}_n \leq \alpha_{ns} \leq y_s \overline{r}^{u}_n, \quad \forall n,s  \label{eq:SAA_manual} \\
& \beta_n, r^{+}_{ns}, r^{-}_{ns} \geq 0, \quad \forall n,s \label{eq:SAA_character}\\
& \sum_{s \in \mathcal{S}} y_s \leq q \label{eq:SAA_violation}\\
%Constraint12
& y_s \in \{0,1\}, \quad \forall s,\label{eq:SAA_bincharacter}
\end{align}
\end{subequations}
\noindent where the set of decision variables is $\Theta = (p_n, r^{-}_{ns}, r^{+}_{ns}$ $r^{d}_{g}, r^{u}_{g}, y_s, \alpha_{ns}, \beta_n)$. Constraints~\eqref{eq:SAA_manual}--\eqref{eq:SAA_bincharacter} represent the sample-based MIP reformulation of the joint chance-constraint \eqref{eq:manual_FJCCOPF}. For a given scenario $s\in \mathcal{S}$, the inequalities \eqref{eq:SAA_manual} establish that a manual adjustment to the production of the dispatchable unit at node $n$ in scenario $s$ can only be done when $y_s = 1$. Otherwise, if $y_s=0$, the power forecast error must be handled by the AGC. Expression \eqref{eq:SAA_violation} ensures that the probability of the joint chance-constraint is met and \eqref{eq:SAA_bincharacter} enforces the binary character of variables $y_s$. A MIP reformulation for the sample average approximation of  \eqref{eq:JCC-OPF} can be found in~\cite{porras2022tight}.

As mentioned above, problem~\eqref{eq:SAA_FJCC-OPF} (and the SAA-based reformulation of \eqref{eq:JCC-OPF} too) becomes rapidly intractable as the power system size and/or the number of scenarios $|\mathcal{S}|$ grows. To alleviate this issue, we make use of the methodology proposed in \cite{porras2022tight}, where a synergistic combination of constraint screening and valid inequalities is exploited. In particular,  line flow constraints that remain nonbinding for all scenarios $s \in \mathcal{S}$ have been screened out. {\color{black} The identification of these constraints is carried out by computing the worst-case flow per line $l$ and per scenario $s$ in both directions over a relaxation of problem \eqref{eq:SAA_FJCC-OPF}. Essentially, this involves solving $2 \times |\mathcal{S}|$ linear programs per line $l$ maximizing (minimizing) the corresponding line flow. If the respective maximum (minimum) line flow in scenario $s$ is strictly smaller (greater) than (minus) the line capacity $\overline{f}_{l}$ ($-\overline{f}_{l}$), the $\leq$-constraint ($\geq$) in~\eqref{eq:SAA_flow} can be safely removed from~\eqref{eq:SAA_FJCC-OPF}.} Furthermore, to tighten the relaxed LP formulation of~\eqref{eq:SAA_FJCC-OPF}, we introduce the valid inequalities developed in \cite{porras2022tight}, which guarantee that, in at least $(1-\epsilon)$-percentage of scenarios, constraints~\eqref{eq:SAA_resca} and~\eqref{eq:SAA_flow} must be individually satisfied by the deployment of AGC only.

Alternatively, we have also implemented an adaptation of the inner convex approximation ALSO-X \cite{jiang2022also}, which is able to identify good feasible solutions of problem~\eqref{eq:SAA_FJCC-OPF}. As the original approximation algorithm, our adaptation works iteratively: It first relaxes the integrality of $\boldsymbol{y}$ and then enters a loop whose core step is to perform a bisection search over parameter $q$. A pseudocode is provided in Algorithm \ref{alg:adapt_alsox}.

\begin{algorithm}
\caption{Adaptation of ALSO-X}\label{alg:adapt_alsox}
Input: Stopping tolerance parameter $\delta$
\begin{algorithmic}[1]
\Require Relax the integrality of $\boldsymbol{y}$
\State $\underline{q} \gets 0$, $\overline{q} \gets \floor*{\epsilon |\mathcal{S}|}$
\While{$\overline{q} - \underline{q} \geq \delta$}
    \State Set $q=(\underline{q}+\overline{q})/2$ and retrieve $\Theta^{*}$ as an optimal solution to \eqref{eq:SAA_FJCC-OPF}.
    \State Set $\underline{q} = q$ if $\mathbb{P}\left(\boldsymbol{y}^{*} = 0 \right)\geq 1-\epsilon$; otherwise, $\overline{q} = q$
\EndWhile
\end{algorithmic}
Output: A feasible solution of model \eqref{eq:SAA_FJCC-OPF}.
\end{algorithm}

%%%%%%%%%%%%%%%%%%%%%%%%%%%%%%%%%%%%%%%%%%%%%%%%
\section{Evaluation  procedure}
\label{sec:eval}
In this section, we outline the procedure for evaluating the performance of the two approaches compared in this paper, namely:
\begin{itemize}
    \item[-] The joint chance-constrained problem with automatic generation control only formulated in \eqref{eq:JCC-OPF} and denoted as AGC-$\epsilon$ hereinafter. For $\epsilon=0$, the constrains must be satisfied for all scenarios and model \eqref{eq:JCC-OPF} is formulated as a linear program. For $\epsilon\neq 0$, model \eqref{eq:JCC-OPF} is reformulated as a MIP problem and efficiently solved using the procedure described in \cite{porras2022tight}.
    \item[-] The joint chance-constrained problem with both automatic and manual generation control formulated in \eqref{eq:FJCC-OPF} and denoted as AMGC-$\epsilon$. Notice that for $\epsilon=0$, the results obtained by AGC-0 and AMGC-0 are the same. For $\epsilon\neq 0$, model \eqref{eq:FJCC-OPF} is reformulated as the MIP model \eqref{eq:SAA_FJCC-OPF} and solved using the procedure described in Section \ref{sec:solution}. The approach that solves model \eqref{eq:SAA_FJCC-OPF} using the heuristic procedure described in Algorithm \ref{alg:adapt_alsox} is denoted as AMGC-H-$\epsilon$.
\end{itemize}

First, let $(p_n^*, r^{d,*}_{g}, r^{u,*}_{g}, \beta^*_n)$ denote the optimal dispatch and reserve capacity decisions delivered by AGC-$\epsilon$, AMGC-$\epsilon$ or AMGC-H-$\epsilon$ with the in-sample scenario set $\mathcal{S}$. We evaluate the performance of these decisions on an out-of-sample scenario set denoted by $\mathcal{S}'$ and indexed by $s'$, with $|\mathcal{S}|\ll|\mathcal{S}'|$. Each out-of-sample scenario $s'$ is characterized by the realization of the forecast errors $\omega_{ns'}$ and the system-wise aggregate forecast error $\Omega_{s'}$, with $\Omega_{s'}=\sum_{n\in\mathcal{N}}\omega_{ns'}$. For each scenario $s'$, we formulate the following real-time operation problem:
\begin{subequations}\label{eq:out_manual}
    \begin{align}
    \min_{\Psi} \hspace{0.1cm} &  \sum_{n} c_{n} \, p^*_{n} + c^{u}_{n} r^{u,*}_{n} + c^{d}_{n} r^{d,*}_{n} + c^{+}_nr^{+}_{ns'} \notag\\
    &   \hspace{2.5cm} - c^{-}_n r^{-}_{ns'} + P(\Delta^{+}_{ns'}  +  \Delta^{-}_{ns'})\label{eq:out_of}\\
    % Power Balance
    \text{s.t.} \hspace{0.1cm} & \sum_{n} \alpha_{ns'} + \Delta^{+}_{ns'} - \Delta^{-}_{ns'} = 0, \label{eq:out_balance}\\
    & r^{+}_{ns'} - r^{-}_{ns'} = -\Omega_{s'}\beta^{*}_n + \alpha_{ns'}, \quad \forall n \label{eq:out_updownreserve}\\
    & -r^{d,*}_{n} \leq  -\Omega_{s'}\beta^{*}_n + \alpha_{ns'} \leq r^{u,*}_{n}, \quad \forall n
     \label{eq:out_resca}\\
    & -\overline{f}_{l} \leq \sum_{n} B_{ln} (p^*_{n} -\Omega_{s'}\beta^*_n + \alpha_{ns'} + \widehat{w}_n +\notag\\
    & \hspace{1cm}  + \omega_{ns'} - d_n + \Delta^{+}_{ns'} - \Delta^{-}_{ns'}) \leq \overline{f}_{l}, \quad \forall l\label{eq:out_flow}\\
    &r^{+}_{ns'},r^{-}_{ns'},\Delta^{+}_{ns'},\Delta^{-}_{ns'} \geq 0, \quad \forall n \label{eq:out_character}.
    \end{align}
\end{subequations}
Note that $\Psi = ( r^{+}_{ns'},r^{-}_{ns'},\alpha_{ns'},\Delta^{+}_{ns'},\Delta^{-}_{ns'})$ where $\Delta^+_{ns'}$ and $\Delta^-_{ns'}$ are two slack variables that quantify the positive and negative power deviations at each node $n$, respectively. These deviations are penalized in the objective function through parameter $P$, {\color{black}which is to be set large enough so that scenario $s'$ is counted as \emph{infeasible} if any of the corresponding slack variables takes on a strictly positive value}. For each scenario $s'$, model \eqref{eq:out_manual} determines the reserve deployments to keep the network balanced at the minimum cost. Note that constraints \eqref{eq:out_balance}-\eqref{eq:out_character} are equivalent to constraints \eqref{eq:SAA_balance}-\eqref{eq:SAA_character} but with the addition of the slack variables $\Delta^{+}_{ns'},\Delta^{-}_{ns'}$, which guarantee the feasibility of model \eqref{eq:out_manual} for any scenario realization.

Using the solution to model \eqref{eq:out_manual}, we split the out-of-sample scenario set $\mathcal{S}'$ into three subsets as follows. First, we solve model \eqref{eq:out_manual} with variables $\alpha_{ns'}$, $\Delta^{+}_{ns'}$ and $\Delta^{-}_{ns'}$ fixed to 0, that is, enforcing that forecast errors can only be handled by AGC. If the problem is feasible, the scenario $s'$ belongs to subset $\mathcal{S}'_{\text A}$ and the real-time operation cost is denoted by $C_{s'}$. If this problem is infeasible, we solve model \eqref{eq:out_manual} without fixing any variable. If $\max(\Delta^+_{ns'},\Delta^-_{ns'})=0$, then the forecast errors can be offset using automatic and manual reserves, and scenario $s'$ is included in subset $\mathcal{S}'_{\text M}$. Finally, if $\max(\Delta^+_{ns'},\Delta^-_{ns'})>0$, then automatic and manual reserves are not enough to maintain the power balance and power deviations occur during the real-time operation. In that case,  scenario $s'$ belongs to the subset $\mathcal{S}'_{\text D}$. In the next section, we evaluate the performance of the different approaches by comparing the percentage of scenarios that belong to each subset as well as the expected cost of the real-time operation.

\section{Numerical Results}
\label{sec:case}

In this section, we compare the performance of  the different approaches presented in Section \ref{sec:formulation} using {\color{black} modified versions of the IEEE-118 and  IEEE-300 test systems widely employed in the technical literature on the topic. The 118-bus system (medium case) shows the performance of our novel model in a worst-case scenario where very few generators can deploy reserve and either have a very cheap or very expensive marginal cost. On the other hand, the 300-bus system (large case) illustrates its performance in a larger system where, in addition, more flexibility is available by the presence of more generators with the ability to back up.}

\subsection{Medium Case: IEEE-118}
The IEEE-118 test system has 118 nodes, 19 generators and 186 transmission lines, and the original data pertaining to this system are publicly available in the repository \cite{pglib}. We assume that six generators can provide reserves, and their corresponding data is shown in Table \ref{tab:generators_reserve}. Notice that units 12, 65 and 111 have a much higher production cost than units 49, 61 and 100. For these six units, the reserve deployment costs are computed as $c^{-}_n = 0.8 \, c_n$ and $c^{+}_n = 1.2 \, c_n$, and the reserve capacity costs are $c^d_n = c^u_n = 0.2c_n$. Besides, we add 25 wind power plants throughout the system as suggested in \cite{roald2016chance}. We consider that the wind power forecast error is normally distributed, i.e., $\boldsymbol{\omega} \sim N(\mathbf{0},\Sigma)$, where $\mathbf{0}$ and $\Sigma$ represent, respectively, the zero vector and the covariance matrix. We also assume that the standard deviation of $\omega_n$ at node $n$ is a 15\% of the wind power forecast $w_n$. The uncertainty pertaining to the renewable generation of the wind farms is characterized using 1000 scenarios, that is, $|\mathcal{S}|=1000$. Finally, the penalty cost $P$ due to deviations is twice the production cost of the most expensive generator. All data of this modified 118-bus system is available at \cite{AGCmanual2023}.

\begin{table}[h]
\begin{center}
\caption{Generators with capability to provide reserve.}
\begin{tabular}{cccccc}
\hline
 $n$ &$c_n$  &$c^{-}_n$ &$c^{+}_n$ &$c^{d}_n/c^{u}_n$ & $\overline{r}^{d}_n/\overline{r}^{u}_n$ \\
\hline
12 &124.6 &99.7	&149.5 &24.9 &85\\
49 &16.7 &13.3	&20.0 &3.3 &223\\
61 &16.1 &12.8	&19.3 &3.2 &195\\
65 &100.0 &80.0	&120.0 &20.0 &441\\
100 &12.6 &10.1	&15.1 &2.5 &653\\
111 &110.0 &88.0	&132.0 &22.0 &79\\
\hline
\end{tabular}
\label{tab:generators_reserve}
\end{center}
\end{table}

%Generators 2, 7, 10, 11, 17 and 19 have the ability to provide reserve, i.e., $\overline{r}^{d}_n$ and $\overline{r}^{u}_n$ are different from 0 for generator $n$ where the value used for each parameter is $\overline{p}_n$. For illustrative purposes, the linear cost of generators 10 and 11 is changed to 100 and 110, respectively; thus, there are 3 cheap generators (with a linear operating cost below 20\euro/MWh) and 3 expensive generators (with a linear operating cost above 100\euro/MWh) providing reserve and, then, a conservative decision implies a high cost. The cost of reserve capacity is a 20\% of the linear operating cost, and we choose $c^{-}_n = 0.8 \, c_n$ and $c^{+}_n = 1.2 \, c_n$ for generator $n$  as the down- and up-reserve deployment cost, in that order.

%, concretely, we compare the performance of the classical JCC-OPF with the FJCC-OPF. To do that, for the JCC-OPF, we use the SAA, ALSO-X and SA reformulation, referred as SAA-JCC-OPF, ALSOX-JCC-OPF and SA-JCC-OPF, in that order (mentioned based on its conservatiness); and, for the FJCC-OPF, we utilize the SAA and MALSO-X reformulation, referred as SAA-JCC-OPF and MALSOX-JCC-OPF, respectively.

To provide meaningful statistics, each method is run for ten different sets of randomly generated scenarios. Accordingly,  we report results averaged over these ten instances. All optimization problems have been solved using GUROBI 9.1.2 \cite{gurobi} on a Linux-based server with CPUs clocking at 2.6 GHz, 6 threads and 16 GB of RAM. In all cases, the optimality GAP has been set to $10^{-9}\%$ and the time limit to 10 hours.

We compare the different approaches following the out-of-sample evaluation procedure described in Section \ref{sec:eval} with 100\,000 different scenarios drawn from the same distribution, that is, $|\mathcal{S}'|=100\,000$. We compare the results of four different approaches, namely: i) AGC-0, which corresponds to the joint chance-constrained OPF model \eqref{eq:JCC-OPF} with $\epsilon=0\%$ (i.e., all scenarios must be satisfied); ii) AGC-5, which is the joint chance-constrained OPF model \eqref{eq:JCC-OPF} with $\epsilon=5\%$ (that is, 50 scenarios may have violated constraints); iii) the proposed AMGC-5 approach, which is the proposed stochastic OPF model \eqref{eq:FJCC-OPF} with $\epsilon=5\%$ (50 scenarios may use manual adjustments to re-dispatch generators), and iv) AMGC-H-5, which is the proposed stochastic OPF model \eqref{eq:FJCC-OPF} solved with the heuristic ALSO-X procedure to solve AMGC-5. The results shown in Table \ref{tab:comparison_scenarios} include (i) the percentage of scenarios in which the forecast errors are handled using automatic generation control only $|\mathcal{S}'_{\text A}|$, (ii) the percentage of scenarios in which manual re-dispatch is required to keep power balance throughout the network $|\mathcal{S}'_{\text M}|$, (iii) the percentage of scenarios in which automatic and manual reserves are not enough to offset power imbalances and therefore, power deviations are inevitable and the system security is compromised $|\mathcal{S}'_{\text D}|$, and (iv) the total expected cost.

%For all numerical simulations presented in this section, the tolerable probability of violation of the joint chance constraint is set to 5\% (i.e., $\epsilon = 0.05$ and $p=50$). This means that, in the case of JCC-OPF, 50 scenarios may have violated constraints, and, in the case of FJCC-OPF, 50 scenarios may make use of manual adjustments to re-dispatch generators.

\begin{table}[h]
\begin{center}
\caption{{\color{black}Medium Case: Out-of-sample performance comparison.}}
\begin{tabular}{lcccccc}
\hline
 &$|\mathcal{S}'_{\text A}|$   &$|\mathcal{S}'_{\text M}|$ &$|\mathcal{S}'_{\text D}|$ & $\mathbb{E}[C_{s'}]$ (k\euro)\\
\hline
\textbf{AGC-0} &99.56\% &0.16\%	&0.28\% &73.65\\
\textbf{AGC-5} &94.12\% &0.17\%	&5.71\% & 53.04\\
\textbf{AMGC-5} &94.30\% &5.41\%	&0.29\% &60.15\\
\textbf{AMGC-H-5} &94.77\% &4.98\%	&0.25\% &61.05\\
\hline
\end{tabular}
\label{tab:comparison_scenarios}
\end{center}
\end{table}

As expected, AGC-0 gives very conservative and expensive solutions, but is able to handle power imbalances using only AGC in 99.56\% of the scenarios and only has 0.28\% scenarios in which system security is compromised. In comparison, we observe that the expected cost of the AGC-5 solutions is 28\% lower on average than those of AGC-0 and that AGC activation is sufficient to guarantee system security in 94.12\% of scenarios (close to the desired violation probability of 95\%). However, the AGC-5 are not able to offset power imbalances with both types of reserves in 5.71\% of the scenarios, thus compromising the safety of the system.
%The ALSOX approach leads to more conservative results than the SAA approach at the expense of increasing the expected cost. Nevertheless, imbalances still occur in 3.71\% of the scenarios.
The proposed AMGC-5 approach is able to compute dispatch solutions while taking into account both automatic and manual reserves. As expected, imbalances are compensated with AGC in 94.30\% of the scenarios, and manual reserves are only required in 5.41\% of them. Note that these percentages are very close to the desired values of 95\% and 5\%. Although there is still 0.29\% scenarios in which system security is compromised, the proposed methodology reduces the expected cost by 18\% with respect to AGC-0 (which has reliability levels). Finally, the results obtained with the heuristic ALSO-X algorithm AMGC-H-5 are slightly more conservative and expensive than those of AMGC-5. However, these results confirm that Algorithm~\ref{alg:adapt_alsox} is able to provide a good feasible solution to the proposed formulation~\eqref{eq:SAA_FJCC-OPF}.

\begin{table}[h]
\begin{center}
\caption{Comparison of decisions made by AGC and AMGC.}
\begin{tabular}{lccccccc}
\hline
&\multicolumn{3}{c}{Cheap Generators} &&\multicolumn{3}{c}{Expensive Generators}\\
 & $\beta$ & $r^{d}$ &$r^{u}$  && $\beta$  & $r^{d}$ & $r^{u}$\\
\hline
\textbf{AGC-0}    &0.77 &533.6 &538.5 &&0.23 &160.1 &162.5\\
\textbf{AGC-5}    &1.00 &509.4 &364.5 &&0.00 & 0.0 &0.0\\
\textbf{AMGC-5}   &0.99 &684.7 &386.3       &&0.01 &7.9 &314.7 \\
\textbf{AMGC-H-5} &0.96 &675.7 &395.3 &&0.04 &16.9 &305.6\\
\hline
\end{tabular}
\label{tab:comparison_decisions}
\end{center}
\end{table}

To give more details on the out-of-sample performance of the different approaches, Table~\ref{tab:comparison_decisions} gathers a summary of the OPF decisions $\beta^*_n$, $r^{d,*}_n$ and $r^{u,*}_n$. For conciseness, we aggregate the units providing reserves into cheap and expensive generators. Interestingly,  AGC-0 yields more conservative OPF decisions since both cheap and expensive generators are dispatched to provide AGC. Conversely, the other methods mainly allocate AGC to cheap generators only. Notice that in the case of AGC-5, the values of $\beta$, $r^d$ and $r^u$ are 0 for expensive generators, which means that these units will not be available for manual reserve during the real-time operation of the power system and therefore, the probability of incurring in dangerous power deviations increases. Conversely, both AMGC-5 and AMGC-H-5 procure more reserve capacities so that cheap and expensive generators can be effectively and efficiently redispatched to minimize the real-time operation cost while reducing power deviations.

\begin{figure}[tbp]
\centerline{\includegraphics[scale=0.4]{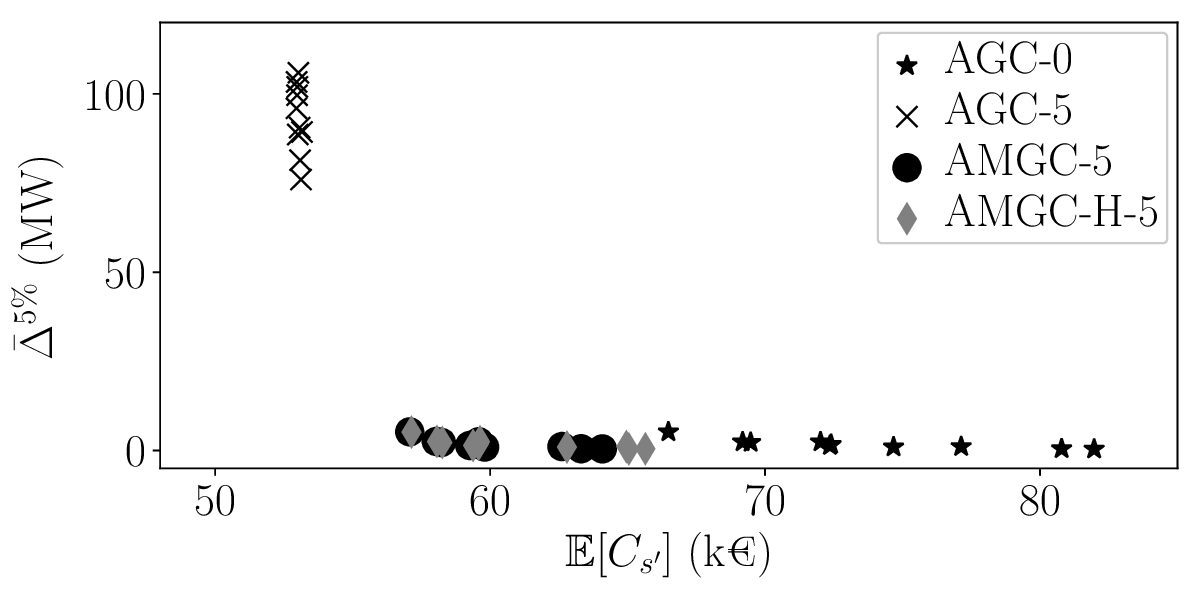}}
\caption{\small Expected cost vs. the average deviation of the \%5 scenarios with the largest deviations \textcolor{black}{for each of the four considered methods. The results for AGC-0 (denoted by stars) are clustered to the lower right of the figure, indicating that the method finds very reliable, but also very costly solutions. The results for AGC-5 (denoted by crosses) are clustered to the upper left of the figure, indicating that this method finds less costly, but also much less reliable solutions.The results for the proposed approach AMGC-5 and the associated heuristic method AMGC-H-5 are clustered in the middle of the figure, indicating that the proposed method finds a good trade-off between cost and reliability. Note that the markers for AMGC-5 and AMGC-H-5 partially overlap, indicating that the heuristic AMGC-H-5 method is able to identify close-to-optimal solutions.}}
\label{fig:Cost_Violation}
\vspace{-0.5cm}
\end{figure}

To further illustrate the differences between the methods compared in this section, we compute for each out-of-sample scenario $s'$ \textcolor{black}{the overall level of infeasibility, quantified by} the total deviation $\Delta_{s'}$ as follows:
\begin{align*}
    & \Delta_{s'} = \sum_n \left(\Delta^{+}_{ns'} + \Delta^{-}_{ns'}\right)
\end{align*}
Figure \ref{fig:Cost_Violation} plots the average value of $\Delta_{s'}$ for the 5\% scenarios with largest deviations ($\bar{\Delta}^{\text 5\%}$) versus the expected cost for each method and each of the 10 independent samples. As observed, AGC-0  involves very low deviation levels but the highest expected cost. Under  AGC-5, the expected cost is decreased at the expense of increasing the system deviations. Finally, the proposed methods AMGC-5 and AMGC-H-5 manage to maintain similarly low levels of deviations as AGC-0, but at a significantly lower cost.

To conclude this case study, we provide the average computational times of the different approaches. Since AGC-0 is formulated as a linear programming problem, it takes 6.1s on average to be solved. Using the efficient solution procedure proposed in \cite{porras2022tight}, AGC-5 is solved in 14.3s. Since the proposed AMGC-5 requires the use of extra variables to properly model the deployment of manual reserves, its average computational time increases up to 3288.5s. Nevertheless, the heuristic procedure described in Section~\ref{sec:solution} is able to reduce this time to 73.6s with a very slight impact on the performance of AMGC-5.

\subsection{Large Case: IEEE-300}
{\color{black} The IEEE-300 test system has 300 nodes, 57 generators, and 411 transmission lines, and the original data pertaining to this system are also publicly available in the repository \cite{pglib}. This system has 15 generators with the ability to use reserve, i.e., it is possible to size the reserve. Small wind farms have been distributed throughout the system, a total of 119, with a share of 40\%. The forecast error is also normally distributed, as in the IEEE-118 case study, and has a standard deviation of 15\% of the predicted value (the optimization setup of the solver, and the number of in-sample and out-of-sample scenarios are the same too). All data of this modified 300-bus system is available at \cite{AGCmanual2023}.}

\begin{table}[h]
\begin{center}
\caption{{\color{black}Large Case: Out-of-sample performance comparison.}}
\begin{tabular}{lcccccc}
\hline
 &$|\mathcal{S}'_{\text A}|$   &$|\mathcal{S}'_{\text M}|$ &$|\mathcal{S}'_{\text D}|$ & $\mathbb{E}[C_{s'}]$ (k\euro)\\
\hline
\textbf{AGC-0} &99.22\% &0.21\%	&0.56\% &242.73\\
\textbf{AGC-5} &93.86\% &0.39\%	&5.75\% & 238.51\\
\textbf{AMGC-5} &94.31\% &5.03\%	&0.54\% &239.64\\
\textbf{AMGC-H-5} &94.78\% &4.68\%	&0.65\% &240.92\\
\hline
\end{tabular}
\label{tab:300_comparison_scenarios}
\end{center}
\end{table}

{\color{black}As can be seen in Table \ref{tab:300_comparison_scenarios}, in this system, results similar to those in Table \ref{tab:comparison_scenarios} are obtained. AGC-0 obtains the most conservative and expensive solution, where 99.22\% of the scenarios are handled by AGC. On the contrary, AGC-5 obtains a solution that is approximately 1.8\% cheaper\footnote{Note that, in this system, the cost difference is much lower because the number of generators with the ability to provide reserve is greater, i.e., it accounts for more flexibility resulting in more optimal affine control policies for AGC-0.} than AGC-0, but is much less reliable as the system's security is not guaranteed in 5.75\% of the scenarios.
%The proposed approach, AMGC-5 performs as expected, and as occurs in the 118 system.
For the proposed approach AMGC-5, imbalances are compensated with AGC in 94.31\% of the scenarios and manual adjustments are required in 5.03\% of them. Its reliability is similar to AGC-0, with only about 0.5\% of samples with system insecurity. However, the AMGC-5 solution is 1.3\% cheaper. This again highlights the ability of the proposed model to provide more reliable and cheaper solutions compared to existing methods.

As mentioned above, an unfortunate aspect of the proposed AMGC-5 model is that it requires significantly higher computational time compared to AGC-0 and AGC-5. Both AGC-0 and AGC-5 can be solved in approximately 20 seconds and 100 seconds, respectively. In contrast, AMGC-5 exceeds the allocated time limit of 10 hours, and terminates with a final MIPGap achieved in the range of 0.2\% to 0.5\%. One option to reduce computation time is to leverage the proposed heuristic method AMGC-H-5 instead of AMGC-5. For the 300-bus system, AMGC-H-5 solves in about 170 seconds. This increase in computation time comes at the expense of slightly more expensive solutions. While AMGC-H-5 is as reliable or more reliable than AMGC-5 and AGC-5, it does have 0.5\% higher operating costs than AMGC-5. Given that this cost is still significantly lower than that of AGC-0 and the reliability of the solution is much better than for AGC-5, we conclude that the heuristic method AMGC-H-5 is a good option for scaling our proposed optimization model to larger systems.}

\section{Conclusions}
\label{sec:conclusion}

Existing approaches to solve the stochastic OPF are either overly conservative and expensive, or leave the system vulnerable to low probability, high impact events. To address this issue, we present a novel stochastic optimal power flow formulation that distinguishes between ``normal'' operation conditions in which power deviations are balanced with AGC only, and ``adverse'' operation under which manual re-dispatch actions are required. As a result, our approach yields solutions that are more reliable and less conservative than existing approaches in the literature.

%In this paper, we propose a novel formulation to address the DC-OPF problem under uncertainty. Our formulation ensures that AGC is able to cover power imbalances, and thus, guarantee system security, with a probability greater than or equal to $1-\epsilon$. For the remaining uncertainty realizations, manual adjustments are modeled to cope with the most adverse scenarios. This is done by means of a joint chance-constraint that limits the probability that operators manually adjust the output dispatch of the generators. As a result, our approach yields solutions that are more reliable and less conservative than existing approaches in the literature.

Our model is formulated as a joint chance-constraint program that limits the probability that operators manually adjust the power output of the generators. {\color{black}To assess the contributions of our proposal, we compare it with existing approaches using an illustrative 3-bus network and more realistic ones such as the 118-bus and 300-bus systems.} The obtained results for the larger system\textcolor{black}{s} demonstrate that the proposed methodology is able to achieve dispatch decisions that maintain almost identical security levels, but are cheaper than approaches that pursue feasibility with AGC actions only under any uncertainty realization. A drawback of our proposed approach is that the computational burden increases due to the modeling of the manual re-dispatch actions. However, we also suggest an heuristic algorithm to solve the proposed model and verify that the computational time is drastically shortened without causing a significant decline in performance.

\textcolor{black}{While the proposed approach represents a step towards more realistic modeling of reserve activation in stochastic OPF problems, several open questions remain. Interesting avenues for future work include for example how to incorporate restrictions on manual reserve activation (such as e.g. limiting the number of generators that participate in the manual activation), a more detailed and realistic modeling of reserve cost that acknowledges the possible cost difference between automatic AGC reserves and manually activated reserves, as well as a strategy to handle situations in which the problem becomes infeasible (due to e.g. high levels of uncertainty, severe system congestion or limited generation capacity).}

%conducted a case study where the performance of the different methods  is evaluated on a test data set. Our results show that the joint chance-constrained Optimal Power Flow, i.e., model \eqref{eq:JCC-OPF} with $\epsilon \neq 0$, exhibits a risky performance, with 5.71\% of the out-of-sample scenarios where the power system is left vulnerable and with its security compromised. On the other hand, the robustly feasible approach, i.e., model \eqref{eq:JCC-OPF} where $\epsilon = 0$, performs adequately guaranteeing system security by means of AGC in 99.56\% of the out-of-sample scenarios. However, this approach relies on a remarkably costly operation. Our approach provides intermediate solutions where AGC and manual reserve are utilized in approximately 95\% and 5\% of the scenarios, respectively. This significantly reduces the operational cost, while maintaining a level of reliability comparable to a fully robust approach. Finally, the proposed heuristic procedure yields similar solutions to the original mixed-integer program, while the associated computational burden makes it competitive with respect to the approaches from the literature.

%%%%%%%%%%%%%%%%%%%%%%%%%%%%%%%%%%%%%%%%%%%%%%%%

\bibliographystyle{ieeetr}
\bibliography{mybibfile}
%%%%%%%%%%%%%%%%%%%%%%%%%%%%%%%%%%%%%%%%%%%%%%%%
\end{document}